\documentclass[11pt, reqno, a4wide]{article}
\usepackage{amssymb, amsmath, amsthm}
\usepackage{ifthen}
\usepackage[all, cmtip]{xy}
\usepackage{pdfpages}
\usepackage{graphics}
\usepackage{mathrsfs}
\usepackage{hyperref}
\newtheorem{leer}{}[section]
\newtheorem{thm}[leer]{Theorem}

\newtheorem{rema}[leer]{Remark} 
\newtheorem{prop}[leer]{Proposition}
\newtheorem{lemm}[leer]{Lemma}
\newtheorem{coro}[leer]{Corollary}
\newtheorem{defi}[leer]{Definition}

\newtheorem*{thm*}{Theorem}
\newtheorem*{conj*}{Conjecture}
\newtheorem*{rema*}{Remark} 
\newtheorem*{prop*}{Proposition}
\newtheorem*{lemm*}{Lemma}
\newtheorem*{coro*}{Corollary}
\newtheorem*{defi*}{Definition}
\newtheorem*{defirema*}{Definition and Remark}
\newtheorem*{exam*}{Example}

\newcommand{\ol}[1]{\overline{#1}}
\newcommand{\ul}[1]{\underline{#1}}

\newcommand{\plim}{\mathop{\varprojlim}\limits}

\newcommand{\EEE}{\mathscr{E}}

\newcommand{\XXX}{\mathscr{X}}

\newcommand{\Cc}{{\mathbb{C}}}

\newcommand{\Ff}{{\mathbb{F}}}
\newcommand{\Gg}{{\mathbb{G}}}

\newcommand{\Ll}{{\mathbb{L}}}

\newcommand{\Nn}{{\mathbb{N}}}

\newcommand{\Qq}{{\mathbb{Q}}}
\newcommand{\Rr}{{\mathbb{R}}}

\newcommand{\Zz}{{\mathbb{Z}}}

\newcommand{\uG}{{\ul{G}}}

\newcommand{\uS}{{\ul{S}}}

\newcommand{\uZ}{{\ul{Z}}}

\newcommand{\oK}{{\ol{K}}}

\newcommand{\GL}{\mathrm{GL}}

\newcommand{\tors}{\mathrm{tors}}

\newcommand{\Gal}{\mathrm{Gal}}
\newcommand{\rk}{\mathrm{rk}}

\newcommand{\ab}{\mathrm{ab}}

\newcommand{\FSQ}{\mathrm{FSQ}}
\newcommand{\JH}{\mathrm{JH}}
\newcommand{\perf}{\mathrm{per}}
\newcommand{\Lie}{\mathrm{Lie}}

\newcommand{\orho}{\overline{\rho}}
\newcommand{\oGamma}{\overline{\Gamma}}
\newcommand{\pr}{\mathrm{pr}}
\renewcommand{\tors}{\mathrm{tor}}
\newcommand{\okappa}{\ol{\kappa}}
\newcommand{\s}{\mathrm{sep}}

\begin{document}
\title{Finiteness properties of torsion fields of abelian varieties}
\author{Wojciech Gajda and Sebastian Petersen}
\maketitle

\begin{abstract} Let $A$ be an abelian variety  defined over a field $K.$ We study finite generation properties of 
the profinite group $\Gal(K_\tors(A)/K)$ and of certain closed normal subgroups thereof, where $K_\tors(A)$ is the torsion field of $A$ over $K$. 
In fact, we establish more general finite generation properties for monodromy groups attached to smooth projective varieties 
via \'etale cohomology.
We apply this in order to give an independent proof and generalizations of a recent result of Checcoli and Dill about small exponent subfields of $K_\tors(A)/K$ in the number field case. 
We also give an application of our finite generation results in the realm of permanence principles for varieties with the weak Hilbert property.
\end{abstract}

\section{Introduction}
For an abelian variety $A$ over a field $K$ we denote by $K(A_{\tors})$ the field obtained from $K$ by adjoining the coordinates of all torsion points in $A(\ol{K})$. 
%For an abelian variety $A$ over a number field $K$ and $e\in\Nn$, Checcoli and Dill (cf. \cite[Theorem 1]{cd}) recently proved that there exists a finite extension $M/K$ with the following property: 
%affirmatively a question of Habegger mentioned in \cite[Section 4]{cd} to the effect whether \cite[Theorem 1]{cd} holds true over
%other base fields $K$, such as
%global function fields or local fields. In fact we shall answer the question affirmatively in the case where $K$ is finitely generated over its prime field\footnote{This
 %includes the case where $K$ is a global function field.}, in the case, where $K$ is finitely generated over an algebraically closed field 
%  and in the case where $K$ is a finitely generated extension of a local field $\kappa$ (the case where $\trdeg(K/\kappa)=0$ and $\chara(\kappa)=0$ of that statement 
 % being done already in \cite[Section 4]{cd}). 
%We think the following finiteness theorems 
%for Galois groups of torsion fields are the right perspective regarding the above question. 
We define $K_\tors(A){:=}K(A_{\tors})\cap K_\s$ where $K_\s$ is the separable closure of $K$. 
For a number field $k,$ we define $k^\dagger=\prod_\ell k_\ell,$ where $\ell$ runs over all prime numbers and 
$k_\ell$  is the compositum of those finite abelian Galois extensions of $k$ that are 
unramified outside $\ell$ and of degree prime to $\ell$. For example, $\Qq^\dagger$ is the compositum of the fields $\Qq(\exp({\frac{2\pi i}{\ell}}))$ for $\ell\in\Ll,$ by Kronecker-Weber theorem.   
The aim of this manuscript is to establish the following theorem about finite generation properties of Galois groups of such torsion fields and to give two applications thereof.

\begin{thm} \label{main-theorem} (cf. Theorem \ref{main-thm} and Remark \ref{expl})
Let $\kappa$ be a field, $K/\kappa$ a finitely generated field extension and $A/K$ an abelian variety.  %Let $K_\tors=K_\tors(A)$. 
Then the following hold true.
\begin{enumerate}
\item[(a)] The profinite group $\Gal(\okappa K_\tors(A)/\okappa K)$ is topologically finitely generated. 
\item[(b)] If  the absolute Galois group $\Gal(\kappa)$ is topologically finitely generated (e.g., when $\kappa$ is a finite or algebraically closed field), then %the profinite group 
$\Gal(K_\tors(A)/ K)$ is topologically finitely generated.
\item[(c)] If  $\kappa$ is a local field, then %the profinite group 
$\Gal(K_\tors(A)/ K)$ is topologically finitely generated.
\item[(d)] If $\kappa$ is a number field, then there exists a finite Galois extension $k/\kappa$ such that %the profinite group
$\Gal(k^\dagger K_\tors(A)/ k^\dagger K)$ is topologically finitely generated.  
\end{enumerate}
\end{thm}

\begin{rema} If in the situation of Theorem \ref{main-theorem} $\kappa$ is a number field, then 
the profinite group $\Gal(K_\tors(A)/K)$ is certainly not topologically finitely generated because $K_\tors(A)$ contains all roots of unity and thus $\Gal(K_\tors(A)/K)$ has an open normal subgroup
of $\hat{\Zz}^\times$ as a quotient.  
\end{rema}

\begin{rema} We will in fact establish more general finite generation properties for monodromy groups attached to smooth projective varieties 
via \'etale cohomology. We refer the reader to Section \ref{sec:et} for the results and do not go into the technical details within the introduction.
\end{rema}

For every field 
extension $\Omega/K$ define $\EEE_e(\Omega/K)$ to be the set of all 
intermetiate fields $F$ of $\Omega/K$ such that $F/K$ is Galois and $\Gal(F/K)$ is a group of exponent $\le e$. 
Our first application of Theorem \ref{main-theorem} adresses a question of Habegger mentioned in Section 4 of the 
recent preprint \cite{cd} of Checcoli and Dill. There Checcoli and Dill estabished the following Theorem 
(cf. \cite[Theorem 1]{cd}):

{\em Let $K$ be a number field and $A/K$ an abelian variety and let $e\in\Nn$. There exists a finite extension $M/K$ such that $F\subset M_{\ab},$ for all $F\in
\EEE_e(K(A_\tors)/K)$.}

Based on Theorem \ref{main-theorem} we give an independent proof and generalization of \cite[Theorem 1]{cd} and results from \cite[Section 4]{cd} as follows.

\begin{coro} %(cf. Theorem \ref{main-thm} and Corollary \ref{tocor}) 
\label{hab1}
Let $\kappa$ be a field. Let $K/\kappa$ be a finitely generated field extension,  $A/K$ an abelian variety and $e\in\Nn$. 
If $\Gal(\kappa)$ is topologically finitely generated (e.g., when $k$ is finite or algebraically closed) or if $\kappa$ is a local field, then there exists a finite separable extension $M/K$ such that $F\subset M,$ for all $F\in\EEE_e(K_\tors(A)/K)$. 
\end{coro}

\begin{coro} %(cf. Theorem \ref{main-thm} and Corollary \ref{tocor}) 
\label{hab2}
Let $K/\Qq$ be a finitely generated field extension, $A/K$ an abelian variety and $e\in\Nn$. 
Then there exists a number field
$k$ and a finite extension $M/K$ such that $F\subset k^\dagger M\subset M_\ab$ for all $F\in \EEE_e(K_\tors(A)/K)$. 
\end{coro}

Our second application of Theorem \ref{main-theorem} concerns permanence principles for varieties that satisfy the weak Hilbert property. 
We shall give a new proof and generalize \cite[Theorem 1.7]{bfp1} considerably. Relations with a conjecture of Zannier \cite[Section 2]{zan} will be 
explained in Remark \ref{rem:zan}. 
We refer the reader to Section \ref{app} for the details.

%\begin{coro} Let $K$ be a number field and $A/K$ an abelian variety. Let $K_\tors=K_\tors(A)$. Then there exists a number field $k$ such that $[Fk^\dagger:Kk^\dagger]<\infty$ for all $F\in \PPP(K_\tors/K)$. \end{coro}

The strategy of proof for Theorem \ref{main-theorem}(a) is as follows. One can construct an adelic Galois representation $\rho: \Gal(K)\to \prod_{\ell\in\Ll} \GL_{2g}(\Zz_\ell)$ such that $G_A(L):=\rho(\Gal(L))$ is isomorphic to
$\Gal(LK_\tors(A)/L)$ for every extension field $L/K$. From \cite{bgp} one gets that some open subgroup $H$ of $G_A(\ol{\kappa} K)$ 
satisfies the following 
technical condition
(+): {\sl For almost all $\ell\in\Ll$ the profinite group $\pr_\ell(H)$ (projection on $\ell$-th factor of the product) 
is generated by its $\ell$-Sylow subgroups.}
In Section \ref{sec:grp} we prove that every closed subgroup of $\prod_{\ell\in\Ll} \GL_{2g}(\Zz_\ell)$ satisfying condition (+) is topologically finitely generated. This then accounts for the proof of parts (a) and (b) of Theorem \ref{main-theorem}, and the proof of Theorem \ref{main-theorem}(d) is similar, relying on \cite{sp} instead of \cite{bgp}. 
The proof of Theorem \ref{main-theorem}(d) in case $|\kappa|<\infty$ is then straightforward. The proof of Theorem \ref{main-theorem}(c) (case where
$\kappa$ is a local field) relies on Theorem \ref{main-theorem}(a) and the potential semistability from \cite{berthelot}. 
The corollaries follow from Theorem \ref{main-theorem} by applying a seminal group theoretical result of Zelmanov \cite{zel} and Wilson \cite{wil}: 
{\sl Every periodic compact (Hausdorff) group
is locally finite.} %Proofs of the corollaries are carried out at the end of Section~\ref{sec:proofs}.

\begin{center}{\bf Acknowledgement}\end{center} 
W.G. thanks Marc Hindry, Paweł Mleczko and participants of SFARA for enlightening discussions on related topics.
S.P. thanks Lior Bary-Soroker, Arno Fehm and Moshe Jarden for interesting discussions and help with some special questions. We thank Arno Fehm in particular for bringing the work \cite{zel} of Zelmanov to our attention. Some of the topics of this paper started in the research of both authors supported by the grant UMO-2018/31/B/ST1/01474 of the National Centre of Sciences of Poland.

\begin{center}{\bf Notation}\end{center}
Let $\Ll$ be the set of all rational prime numbers. For a field $K$ let $\ol{K}$ 
be an algebraic closure of $K$, $K_\s$ (resp. $K_\perf$) the separable (resp. perfect) closure of $K$ inside $\ol{K}$,
and $\Gal(K)=\Gal(K_\s/K)$ the absolute Galois group of $K$. 
 We denote by $K_\ab$ the maximal abelian extension of $K$ in $K_\s$.  A $K$-variety  is a separated algebraic $K$-scheme which is geometrically reduced and 
 geometrically irreducible.
 
For a profinite group $G$ and $\ell\in\Ll$ we let $S^{(\ell)}(G)$ be the normal subgroup topologically generated by 
the  $\ell$-Sylow subgroups of $G$; if $\ell$ is clear from the context we simply write $G^+$ instead of $S^{(\ell)}(G)$, following \cite{bible}. 
We define $\exp(G):=\inf\{n\in\Nn: \mbox{$g^n=1$ for all $g\in G$}\}\in\Nn\cup\{ \infty\}$ to be the exponent of $G$. 
We let $\FSQ(G)$ (resp. $\JH(G)$) be the class of all finite simple quotients of $G$ (resp. of all Jordan-Hölder factors of $G$). 
We let $\Lie_\ell$ be the class of all finite simple groups of Lie type in characteristic $\ell$. 
We  let $\GL_n^{(1)}(\Zz_\ell)$ be the kernel of the natural surjection $\GL_n(\Zz_\ell)\to \GL_n(\Ff_\ell)$.

%The torsion (resp $n$-torsion) of a $\Zz$-module $M$ is denoted by $M_\tors$ (resp. $M[n]$). We put $M[n^\infty]=\bigcup_j M[n^j]$. 
%We furthermore will consider the group schemes $A[n]=\ker(n_A: A\to A)$

\section{Concepts from group theory}\label{sec:grp}
We recall some information about subgroups of $\GL_n(\Ff_\ell)$. Of central importance is the following theorem of Larsen and Pink. We do not state it in its most general form.

\begin{thm}(Larsen and Pink, cf. \cite[Theorem 0.2]{lp}) \label{lpthm} Let $n\in\Nn$. There exists a constant $J'(n)$, depending only on $n$, such that for every $\ell\in\Ll$ 
and every subgroup $\oGamma$of $\GL_n({\ol{\Ff}}_\ell)$
 there are normal subgroups
$\oGamma\triangleright \oGamma_1\triangleright \oGamma_2 \triangleright \oGamma_3$
of $\oGamma$ such that 
\begin{itemize}
\item[(a)] $|\oGamma/\oGamma_1|\le J'(n)$, 
\item[(b)] $\oGamma_1/\oGamma_2=L_1\times\cdots\times L_s$ is a finite product of groups $L_j\in\Lie_\ell$,
\item[(c)] $\oGamma_2/\oGamma_3$ is abelian of order prime to $\ell$ and
\item[(d)] $\oGamma_3$ is an $\ell$-group.
\end{itemize}
\end{thm}

\begin{rema} \label{lprem} The proof of \cite[Theorem 0.2]{lp} in \cite[p. 1155--1156]{lp} actually gives more information. 
Let $\Ff:=\ol{\Ff}_\ell$ and $\uG/\Ff$ the algebraic subgroup of $\GL_{n, \Ff}$ introduced in that proof. Let $\uZ$ be the center of the reductive group 
$\uG_{\mathrm{red}}^\circ=\uG^\circ /\mathrm{Rad}_u(\uG^\circ)$ and $\uS=\uG_{\mathrm{red}}^\circ/\uZ$. 
\begin{enumerate}
\item[(a)] The number $s$ in Theorem \ref{lpthm}(b) satisfies $s\le \dim(S)\le \dim(\GL_n){=}n^2$. 
\item[(b)] The group $\oGamma_2/\oGamma_3$ in Theorem \ref{lpthm}(c) is contained in the torus $\uZ(\Ff)$. Furthermore, as $\Ff$ is algebraically closed, 
$\uZ\cong \Gg_m^h$ for some $h\in\Nn$ with $h\le \dim(\GL_n)=n^2$. 
\end{enumerate} 
\end{rema}

\begin{defi} For a profinite group $\Gamma$ we define $d(\Gamma)$ to be the minimal number $d$ such that $\Gamma$ can be topologically generated by $d$ elements, i.e., 
such that $\Gamma$ contains a dense subgroup that can be generated by $d$ elements.
\end{defi}

A {\em lattice} in $\Qq_\ell^n$ is a free $\Zz_\ell$-submodule $\Lambda$ of $\Qq_\ell^n$ such that the canonical map $\Lambda\otimes_{\Zz_\ell} \Qq_\ell\to \Qq_\ell^n$ is an isomorphism. We identify the group of $\Zz_\ell$-automorphims $\GL_\Lambda(\Zz_\ell)$ of $\Lambda$ with the subgroup 
$\{f\in \GL_n(\Qq_\ell): f(\Lambda)=\Lambda\}$ of $\GL_n(\Qq_\ell)$ in the sequel. 

\begin{rema}\label{lattice}
Let $\Gamma$ be a compact subgroup of $\GL_n(\Qq_\ell)$. Then there exists a lattice $\Lambda$ in $\Qq^n$ such that $\Gamma\subset \GL_\Lambda(\Zz_\ell)$. In 
particular $\Gamma$ is isomorphic to a closed subgroup of $\GL_n(\Zz_\ell)$. 
\end{rema}

\begin{lemm}\label{lpfg} For every $n\in\Nn$ there exists a bound $b(n)$, depending only on $n$,  such that for every $\ell\in\Ll$, every compact subgroup $\Gamma$ of $\GL_n(\Qq_\ell)$ 
is topologically finitely generated with $d(\Gamma)\le b(n)$. 
\end{lemm}

\begin{proof} By Remark \ref{lattice} we can assume that $\Gamma\subset \GL_n(\Zz_\ell)$.
The profinite group $P:=\Gamma\cap \GL_n^{(1)}(\Zz_\ell)$ is topologically finitely generated with $d(P)\le n^2$
(cf. \cite[Prop. 8.1, Ex. 6.3]{klopsch}, \cite{an}). The group $\ol{\Gamma}:=\Gamma/P$ is isomorphic to a subgroup of $\GL_n(\Ff_\ell)$, hence 
Theorem \ref{lpthm} applies to it. 
Let $\oGamma_1, \oGamma_2$ and $\oGamma_3$ be as in that theorem. 
%It is enough to bound $d(\ol{\Gamma})$ only in terms of $n$. 
%Let $J'(n)$ be the constant from \cite[Theorem 0.2]{lp}.  
Clearly $d(\oGamma/\oGamma_1)\le J'(n)$. The groups $L_j$ from Theorem \ref{lpthm}(b) can be generated by two elements
(cf. \cite{steinberg-gen}, \cite[§1]{mart}), and $s\le n^2$ by Remark \ref{lprem}(a).  
Hence $d(\oGamma_1/\oGamma_2)\le 2n^2$. From Remark \ref{lprem}(b)  
$\oGamma_2/\oGamma_3$ is a subgroup of $(\Ff^\times)^h$. The Pontryagin dual of $\Ff^\times$ is pro-cyclic and $\oGamma_2/\oGamma_3$ is a quotient of the Pontryagin dual of 
$(\Ff^\times)^h$, hence $d(\oGamma_2/\oGamma_3)\le h\le n^2$. 
 Finally 
$d(\oGamma_3)\le \frac{1}{4}n^2$ by \cite[Theorem B]{ap}. 
\end{proof}

\begin{coro} \label{fsq} Let $\ell\in\Ll$ and let $\Gamma$ be a compact subgroup of $\GL_n(\Qq_\ell)$. If $\ell>J'(n)$ and $\Gamma=\Gamma^+$, then $\FSQ(\Gamma)\subset 
\{\Zz/\ell\}\cup \Lie_\ell$. 
\end{coro}

\begin{proof} By Remark \ref{lattice} we can assume that $\Gamma\subset \GL_n(\Zz_\ell)$. The group  $P:=\Gamma\cap \GL_n^{(1)}(\Zz_\ell)$ is pro-$\ell$ and $\oGamma=\Gamma/P$ is isomorphic to a subgroup of $\GL_n(\Ff_\ell)$ satisfying
$\oGamma=\oGamma^+$. Hence the assertion is immediate from Theorem \ref{lpthm}. 
\end{proof}

%For later use let us also recall a theorem of E. Artin. In the situation of Theorem \ref{lpthm}, if $\oGamma=\oGamma^+$ and $\ell>J'(n)$, then $\FSQ(\oGamma)\s

\begin{thm} (cf. \cite{artin}, \cite{KLST}, \cite[Th\'eorème 5]{bible})
\label{artin} If $5\le \ell_1<\ell_2$, then $\Lie_{\ell_1}\cap \Lie_{\ell_2}=\emptyset$.
\end{thm}

%The following definition is inspired by \cite{bible} and \cite{sp}.

\begin{defi}  Let $n\in\Nn$ and $L\subset \Ll$. Let $G$ be a compact subgroup of $\prod_{\ell\in L} \GL_n(\Qq_\ell)$ and $pr_\ell$ the projection on the $\ell$-th factor. 
\begin{enumerate}
\item[(a)] We call $G$ independent (resp. group theoretically independent) if $G=\prod_{\ell\in L} \pr_\ell(G)$ 
(resp. if\;  $\FSQ(\pr_{\ell_1}(G))\cap \FSQ(\pr_{\ell_2}(G))=\emptyset$ for all $\ell_1\neq \ell_2$ in $L$). 
%We use ``gt-independent'' as a shortcut for ``group-theoretically independent''. 
\item[(b)] We say that $G$ satisfies condition (+) if 
$\pr_\ell(G)=\pr_\ell(G)^+$ for all but finitely many $\ell$ in $L$. 
\item[(c)] We say that $G$ satisfies condition (+) potentially if $G$ has an open subgroup $H$ such that $H$ satisfies condition (+)
\end{enumerate}
\end{defi}

\begin{rema} If $G$ is group theoretically independent, then $G$ is independent (cf. \cite[Lemme 2]{bible}). \label{lem3-serre}
\end{rema}

\begin{lemm}  \label{gtfg} Every group theoretically independent compact subgroup $G$ of $\prod_{\ell\in L} \GL_n(\Qq_\ell)$ is topologically finitely generated. 
\end{lemm}

\begin{proof} Let $b=b(n)$ be the constant from Lemma \ref{lpfg}. For every $\ell\in L$ there exists a system  $(g^{(\ell)}_1,\cdots, g^{(\ell)}_b)$ of topological generators of $\pr_\ell(G)$
by Lemma \ref{lpfg}. Consider the $g_j=(g^{(\ell)}_j)_{\ell\in L}\in \prod_{\ell\in L} \pr_\ell(G)$ and the closure $H$ of the subgroup $\langle g_1,\cdots, g_b\rangle$ generated
by the elements $g_j$. Then $\pr_\ell(G)=\pr_\ell(H)$ for all $\ell\in \Ll$. From this and our assumption on $G$ we see that the groups $G$ and $H$ are both group theoretically independent. 
It follows that $G$ and $H$ are independent (cf. Remark \ref{lem3-serre}), and thus $G=H$, as desired. 
\end{proof}

\begin{lemm} \label{kickls} Let $G$ be a compact subgroup of $\prod_{\ell\in L} \GL_n(\Qq_\ell)$, $L_0$ a finite subset of $L$ and 
$\pr: \prod_{\ell\in L} \GL_n(\Qq_\ell)\to \prod_{\ell\in L\setminus L_0} \GL_n(\Qq_\ell)$ the projection. If $\pr(G)$ is topologically finitely generated, then
$G$ is topologically finitely generated. 
\end{lemm}

\begin{proof} $\ker(\pr)\cap G$ is isomorphic to a compact subgroup of the {\em finite} product 
$\prod_{\ell\in L_0} \GL_n(\Qq_\ell)$. It follows from Lemma \ref{lpfg} that $\ker(\pr)\cap G$ is topologically finitely generated. The assertion is now immediate from the exact sequence $1\to \ker(\pr)\cap G\to G\to \pr(G)\to 1$. 
\end{proof} 

\begin{prop} \label{plus} If a compact subgroup $G$ of $\prod_{\ell\in L} \GL_n(\Qq_\ell)$ satisfies condition (+), then it is topologically finitely generated.
\end{prop}

\begin{proof} Let $G_\ell=\pr_\ell(G)$. 
There exists a finite subset $L_0$ of $L$ such that $G_\ell=G_\ell^+$ for all $\ell\in L\setminus L_0$. We can furthermore assume that $L_0$ contains
all rational primes $\le J'(n)$ and the primes $2$ and $3$. For all $\ell\in L\setminus L_0$ we have $\FSQ(G_\ell)=\Lie_\ell\cup \{\Zz/\ell\}$ by Corollary \ref{fsq}. 
By Theorem \ref{artin} we see 
that $\FSQ(G_{\ell_1})\cap \FSQ(G_{\ell_2})=\emptyset$ for all $\ell_1,\ell_2 \in L\setminus L_0$ with $\ell_1\neq \ell_2$. Hence the image $\pr(G)$ of $G$ under the projection
$\pr: \prod_{\ell\in L} \GL_n(\Qq_\ell)\to \prod_{\ell\in L\setminus L_0} \GL_n(\Qq_\ell)$
is group theoretically independent. By Lemma \ref{gtfg} the profinite group $\pr(G)$ is topologically finitely generated, and this suffices by Lemma \ref{kickls}
\end{proof}

For further use, we finally discuss in which circumstances propery (+) descends to normal subgroups. 

\begin{defi} For a profinite group $G$ define $\Ll^{\Lie}(G)$ to be the set of all $\ell\in\Ll$ such that $\JH(G)\cap \Lie_\ell\neq \emptyset$. 
\end{defi}

We note that $\Ll^\Lie(G)$ is finite for example when $|G|<\infty$ or when $G$ is pro-solvable. 

\begin{lemm} \label{plussolv} Let $G$  be a compact subgroup of  $\prod_{\ell\in L} \GL_n(\Qq_\ell)$ and $N$ a closed normal subgroup of $G$. 
\begin{enumerate}
\item[(a)] If $G$ satisfies condition (+) and $\Ll^\Lie(G/N)$ is finite, then $N$ satisfies condition (+).
\item[(b)] If $G$ satisfies condition (+) potentially, then there exists an open {\em normal} subgroup $H$ of $G$ such that
$H$ satisfies condition (+). 
\item[(c)] If $G$ satisfies condition (+) potentially and $\Ll^\Lie(G/N)$ is finite, then $N$ satisfies condition (+) potentially and is topologically finitely generated.
\end{enumerate}
\end{lemm} 

\begin{proof} Assume througout that $\Ll^\Lie(G/N)$ is finite. Let $G_\ell=\pr_\ell(G)$ and $N_\ell=\pr_\ell(N)$. Assume $G$ satisfies condition (+). Then, for all but finitely many $\ell$, we
have $G_\ell=G^+_\ell$ and $\JH(G_\ell/N_\ell)\cap \Lie_\ell=\emptyset$, so that \cite[Lemma 1.6]{sp} implies $N_\ell=N_\ell^+$, whence $N$ satisfies condition (+). This proves (a).  
Now assume that $G$ satisfies condition (+) potentially. Then there exists an open subgroup $H_0$ of $G$ such that $H_0$ satisfies condition (+). If we let 
$H=\bigcap_{g\in G} g^{-1} H_0 g$, then $H$ is an open normal subgroup of $G$ and satisfies (+) by (a). Thus (b) holds true.
Furthermore $H/N\cap H$ is a normal subgroup of $G/N$. It thus follows that $\Ll^\Lie(H/N\cap H)$ is finite, hence  (a) implies that that $N\cap H$ satisfies condition (+). 
As
$N\cap H$ is open in $N$, it follows that $N$ satisfies condition (+) potentially. Lemma \ref{plus} implies that $N\cap H$ is topologically finitely generated. As $N\cap H$ is open in $N$, it follows that $N$ is topologically finitely generated. 
This finishes up the proof of (c). 
\end{proof}

\section{Representations attached to cohomology}\label{sec:et}
Throughout this section $\kappa$ is a field of characteristic $p\ge 0$, $K/\kappa$ a finitely generated extension and $\Ll'=\Ll\setminus \{p\}$. Let
$X/K$ be a smooth projective variety. Let $i\in \Nn$, 
$j\in \Zz$. For every $\ell\in \Ll'$ consider the $\ell$-adic \'etale cohomology group $V_\ell=H^i(X_{\ol{K}}, \Qq_\ell(j))$.  Consider the representaion $\rho_\ell: \Gal(K)\to \GL_{V_\ell}(\Qq_\ell)$.
Let $\rho: \Gal(K)\to \prod_{\ell\in\Ll'} \GL_{V_\ell}(\Qq_\ell)$ be the homorphism induced by the $\rho_\ell$. For every field extension $E/K$ there is 
a restriction map $r_{E/K}: \Gal(E)\to \Gal(K)$ and we define $G(E)=\rho(r_{E/K}(\Gal(E)))$. 

\begin{lemm} \label{residual} Let $\ell\in\Ll'$ and let $E/K$ be a separable algebraic field extension. Let $\varepsilon_\ell: \Gal(K)\to \Qq_\ell^\times$ be the cyclotomic character and let 
$\rho'_\ell=\rho_\ell\otimes \varepsilon_\ell: \Gal(K)\to \GL_{V_\ell}(\Qq_\ell)$. 
If $E$ contains the $\ell$-th roots of unity and $\rho_\ell(\Gal(E))=\rho_\ell(\Gal(E))^+$, then  $\rho'_\ell(\Gal(E))=\rho'_\ell(\Gal(E))^+$.
\end{lemm}

\begin{proof} By Lemma \ref{lattice} there exists a lattice $T$ in $V_\ell$ such that 
$\rho_\ell(\Gal(K))\subset \GL_T(\Zz_\ell)$. From $\varepsilon_\ell(\Gal(K))\subset \Zz_\ell^\times$ we conclude that 
$\rho_\ell'(\Gal(K))\subset \GL_T(\Zz_\ell)$. Let $p: \GL_T(\Zz_\ell)\to \GL_T(\Ff_\ell)$ be the projection and consider the residual representations 
$\orho_\ell=p\circ \rho_\ell$ and $\orho'_\ell=p\circ \rho'_\ell$. From $\rho_\ell(\Gal(E))=\rho_\ell(\Gal(E))^+$ it follows that 
$\orho_\ell(\Gal(E))=\orho_\ell(\Gal(E))^+$. This implies $\orho'_\ell(\Gal(E))=\orho'_\ell(\Gal(E))^+$ because $\varepsilon_\ell(\Gal(E))\subset 1+\ell\Zz_\ell$ by our assumption
on $E$. From this it follows that $\rho'_\ell(\Gal(E))=\rho'_\ell(\Gal(E))^+$ because $\ker(p)$ is pro-$\ell$. 
\end{proof}

There exists $n\in\Nn$ such that $\dim_{\Qq_\ell}(V_\ell)=n$ for all $\ell\in \Ll'$ by the Weil conjectures (cf. \cite[Thm. 1.6]{deligne}, \cite[Rem. 1.4]{illusie}). We can thus 
choose isomorphisms $\GL_{V_\ell}(\Qq_\ell)\cong \GL_n(\Qq_\ell)$ and apply results from Section \ref{sec:grp}. 

\begin{prop} (cf. \cite[Theorem 7.5]{bgp}, \cite[Theorem 3.1]{sp}) \label{p1-et} 
\begin{enumerate}
\item[(a)] The profinite group $G(\okappa K)$  satisfies condition (+) potentially. 
\item[(b)] If $\kappa$ is a number field, then there exists a finite Galois extension $k/\kappa$ such that $G(k^\dagger K)$ satisfies condition (+) potentially. 
\end{enumerate}
\end{prop}

\begin{proof} Part (a)  in case $j=0$ is immediate from \cite[Theorem 7.5]{bgp}.  As $\okappa$ contains all roots of unity $H^i(X_{\oK}, \Qq_\ell)$ and
$H^i(X_{\oK}, \Qq_\ell(j))$ are isomorphic as $\Gal(\okappa K)$-modules. Thus part (a) follows in general. 

From now on assume that $\kappa$ is a number field. Part (b) in case $j=0$ is established in  \cite[Theorem 3.1]{sp}. By Lemma \ref{residual} part (b) follows in general.
\end{proof}

\begin{thm} \label{main-thm-et} The profinite group $G(\okappa K)$ is 
topologically finitely generated. If $\Gal(\kappa)$ is topologically finitely generated (e.g., when $\kappa$ is finite or algebraically closed), then $G(K)$ is 
topologically finitely generated.
\end{thm}

\begin{proof} By Proposition \ref{p1-et}(a) and Lemma \ref{plussolv}(c) the profinite group $G(\okappa K)$ is 
topologically finitely generated. There exists an epimorphism 
$$\Gal(\kappa_\s K/K)\to G(K)/G(\okappa K)$$ 
and the profinite group
$\Gal(\kappa_\s K/K)$ is isomorphic to an
open subgroup of $\Gal(\kappa)$. Hence, if $\Gal(\kappa)$ is topologically finitely generated, then $G(K)/G(\okappa K)$ is topologically finitely generated and it follows that $G(K)$ is finitely generated. 
\end{proof} 

\begin{lemm}\label{lemm-local} If $\kappa$ is a local field and $K=\kappa$, then $G(K)$ is topologically finitely generated.
\end{lemm}

\begin{proof} If $K\in\{\Rr, \Cc\}$ is an archimedian local field, then $\Gal(K)$ is finite and thus the assertion is trivially satisfied by Theorem
\ref{main-thm-et}. So assume that $K$ is a non-archimedian local field.  Let $q$ be the residue characteristic of the local field $K$.
Let $L=\Ll'\setminus \{q\}$, $\rho^*: \Gal(K)\to \prod_{\ell\in L} \GL_{V_\ell}(\Qq_\ell)$ the homomorphism induced by the $\rho_\ell$ for $\ell\in L$ and $G^*(K)=\rho^*(\Gal(K))$. By Lemma
\ref{kickls} it is enough to show that $G^*(K)$ is topologically finitely generated.
Let $I\subset \Gal(K)$ be the inertia group and $P$ the maximal normal pro-$q$ subgroup of $I$. 
 By the semistable reduction theorem \cite[Prop. 6.3.2]{berthelot} there exists an open subgroup $J$ of $I$ 
such that for every $\ell\in L$ the action of $J$ on $H^i(X_{\ol{K}}, \Qq_\ell)$ is unipotent. Furthermore $H^i(X_{\ol{K}}, \Qq_\ell)$ is isomorphic to 
$H^i(X_{\ol{K}}, \Qq_\ell(j))$ as a $J$-module because $\varepsilon_\ell(J)=\{1\}$. It follows that the action of $J$ on $H^i(X_{\ol{K}}, \Qq_\ell(j))$ is unipotent. Thus
$\rho_\ell(J\cap P)=0$ for all $\ell\in L$. Hence $\rho^*(P)$ is finite. As $\Gal(K)/I$ and $I/P$ are topologically finitely generated, it follows that $G^*(K)$ is 
topologically finitely generated, as desired.
\end{proof}

\begin{thm} \label{main-thm-et2} If $\kappa$ is a local field, then $G(K)$ is topologically finitely generated.
\end{thm}

\begin{proof} 
After replacing $\kappa$ by a finite extension (and replacing the rest accordingly) we can assume that $K/\kappa$ is separable 
(cf. \cite[4.6.7]{EGAIV2}) and primary. 
Then there exists a geometrically connected smooth $\kappa$-scheme $S$ with function field $K$. By the usual spreading-out principles, after replacing $S$ by a non-empty open subscheme, we can assume that $X$ extends to a smooth projective $S$-scheme $\XXX$ such that $f: \XXX\to S$ has geometrically connected fibres and such that for every $\ell\in\Ll'$ the sheaf
$R^if_*\Zz_\ell(j)$ is lisse and of formation compatible with any base change $S'\to S$ (cf. \cite[Cor. 2.6]{illusie}). In particular $\rho_\ell$ factors through $\pi_1(S)$. 
After replacing $\kappa$ by a finite separable extension and replacing the rest
accordingly we can assume that there
exists a point $s\in S(\kappa)$. 
Let $\sigma$ be the section of $\pi_1(S)\to \Gal(\kappa)$ induced by $s$ (well-defined up to conjugation). There is a diagram 
$$\xymatrix{ &&G(K)\\
1\ar[r] & \pi_1(S_{\ol{\kappa}})\ar[r]& \pi_1(S)\ar@/^/[r]\ar[u]_{\rho}  & \Gal(\kappa)\ar[r]\ar@/^/[l]^\sigma & 1
}$$
with exact row. Now $\rho(\pi_1(S_{\ol{\kappa}}))= G(\ol{\kappa} K)$ is topologically finitely generated by Theorem \ref{main-thm-et}. 
By the base change compatibility the representation $\rho_\ell\circ \sigma$ of $\Gal(\kappa)$ is isomorphic to the representation of $\Gal(\kappa)$ on
$H^q(X_{s, \okappa}, \Qq_\ell(j))$ where $X_s=f^{-1}(s)$. Hence $\rho(\sigma(\Gal(\kappa)))$ is 
 topologically finitely generated by Lemma \ref{lemm-local}. From the diagram we see that 
 $G(K)=\rho(\pi_1(S_{\ol{\kappa}}))\cdot \rho(\sigma(\Gal(\kappa)))$, hence $G(K)$ is  topologically finitely generated. 
\end{proof}

\begin{rema} If $\kappa$ is a finite extension of $\Qq_p$, then $\Gal(\kappa)$ is known to be finitely generated (cf. \cite[Theorem 3.4]{js}). Hence Theorem 
\ref{main-thm-et2} in that case follows directly from Theorem \ref{main-thm-et}. If $\kappa$ is a finite extension of $\Ff_p((t))$, however, then $\Gal(\kappa)$ is not finitely 
generated, and thus the arguments from the proof of Theorem \ref{main-thm-et2} are needed essentially in that case. 
\end{rema}

\begin{thm} \label{main-thm-et3} 
If $\kappa$ is a number field, then there exists a finite Galois extension $k/\kappa$ with the following property: If $\Omega/k^\dagger K$ is a Galois extension
and $\Ll^\Lie(\Gal(\Omega/k^\dagger K))$ is finite, then
$G(\Omega)$ is topologically finitely generated.
\end{thm}

\begin{proof} Immediate from Proposition \ref{p1-et}(b) and Lemma \ref{plussolv}(c). 
\end{proof}

\section{Application to torsion fields of abelian varieties}\label{sec:proofs}
Let $K$ be a field of characteristic $p\ge 0$, $\Ll'=\Ll\setminus \{p\}$, $A/K$ an abelian variety and $g=\dim(A)$.
For all $\ell \in\Ll$ (including the case $\ell=p$) we consider the Tate module $T_\ell:=\plim_{j\in\Nn} A[\ell^j](\ol{K})$
of the Barsotti Tate group $A[\ell^\infty]$, put $V_\ell=T_\ell\otimes_{\Zz_\ell} \Qq_\ell$ and note that $T_\ell$ is a free $\Zz_\ell$-module with
$$\rk_{\Zz_\ell}(T_\ell)=\dim_{\Ff_\ell}(A[\ell](\ol{K}))=\left\{\begin{array}{ll} 2g & \mbox{if $\ell\neq p$,}\\
\le 2g & \mbox{if $\ell=p$.}
\end{array}\right.$$
There is a natural action of $\Gal(K_\perf)$ \footnote{For $\ell\neq p$ the passage to $K_\perf$ is not necessary because then
the finite group schemes $A[\ell^j]$ are \'etale over $K$ and  $A[\ell^j](\ol{K})= A[\ell^j](K_\s)$.} on $V_\ell$ and a restriction isomorphism $r_{K_\perf/K}: \Gal(K_\perf)\to \Gal(K)$, so we get a 
representaion $\rho_\ell: \Gal(K)\cong \Gal(K_\perf)\to \GL_{V_\ell}(\Qq_\ell)$.
We consider the homomorphisms $$\rho: \Gal(K)\to \prod_{\ell\in \Ll} \GL_{V_\ell}(\Qq_\ell)\ \mbox{and}\ 
\rho^*: \Gal(K)\to \prod_{\ell\in \Ll'} \GL_{V_\ell}(\Qq_\ell) $$
induced by the $\rho_{\ell}$. For every field extension $E/K$ there is a restriction map $r_{E/K}: \Gal(E)\to \Gal(K)$ and 
we define $G_A(E)=\rho(r_{E/K}(\Gal(E)))$ and $G_A^*(E)=\rho^*(r_{E/K}(\Gal(E)))$. 

\begin{rema} \label{expl} Let $K_\tors(A):=K(A_\tors)\cap K_\s$ and $E/K$ a field extension. The homomorphism 
$\rho$ induces an isomorphism $G_A(E)\cong \Gal(EK_\tors(A)/E)$. 
\end{rema}

\begin{thm} \label{main-thm} Let $\kappa$ be a field, $K/\kappa$ a finitely generated extension and $A/K$ an abelian variety. 
\begin{enumerate}
\item[(a)] The profinite group $G_A(\okappa K)$ is  topologically finitely generated. 
\item[(b)] If $\Gal(\kappa)$ is topologically finitely generated (e.g., when $\kappa$ is algebraically closed or finite), then $G_A(K)$ is topologically finitely generated.
\item[(c)] If  $\kappa$ is a local field, then $G_A(K)$ is topologically finitely generated.
\item[(d)] If $\kappa$ is a number field, then there exists a finite Galois extension $k/\kappa$ with the following property: If $\Omega/k^\dagger K$ is a Galois extension
and $\Ll^\Lie(\Gal(\Omega/k^\dagger K))$ is finite, then
$G_A(\Omega)$ is topologically finitely generated.
\end{enumerate}
\end{thm}

\begin{proof} For every $\ell\in\Ll\setminus \{p\}$ there is an $\Gal(K)$-equivariant isomorphism $V_\ell(A)\cong H^1(A^\vee_{\ol{K}}, \Qq_\ell(1))$, where
$A^\vee$ is the dual abelian variety. Hence the statements (a)-(d) with 
$G_A$ replaced by $G_A^*$ are a consequence of Theorems \ref{main-thm-et},  \ref{main-thm-et2} and  \ref{main-thm-et3}. Then Lemma \ref{kickls} implies that the statements (a)-(d) hold true as they stand. 
\end{proof}

\begin{rema} \label{explO} In the situation of Theorem \ref{main-thm}(d) it follows that $G_A(\Omega)$ is topologically finitely generated for 
$\Omega\in\{k^\dagger K, k_\ab K, (kK)_\ab\}$. We introduced the additional field $\Omega$ in part (d) because we need the flexibility in Section {\ref{app}}. 
\end{rema} 

\begin{proof}[Proof of Theorem \ref{main-theorem}] 
This is an immediate consequence of Remark \ref{expl} and Theorem \ref{main-thm}.  
\end{proof}

\begin{prop} \label{tocor} Let $A/K$ be an abelian variety over a field $K$ and $E/K$ an algebraic extension. Let $e\in\Nn$. Assume that $G_A(E)$ is topologically finitely generated. 
There exists a finite separable extension $M/K$ such that $F\subset EM$ for all $F\in\EEE_e(K_\tors/K)$. 
\end{prop} 

\begin{proof} Let $F_{\max}^{(e)}$ be the compositum of all fields in $\EEE_e(K_\tors/K)$. Then $F_{\max}^{(e)}/K$ is Galois and 
$\Gal(F_{\max}^{(e)}/K)$ is periodic. $\Gal(EF_{\max}^{(e)}/E)$ is a quotient of $G_A(E)$ and isomorphic to a subgroup of $\Gal(F_{\max}^{(e)}/K)$.  Hence the profinite group
$\Gal(EF_{\max}^{(e)}/E)$ is topologically finitely generated and periodic. By 
\cite[Theorem 1]{zel} and \cite[Corollary, p. 58]{wil} the group $\Gal(EF_{\max}^{(e)}/E)$ is finite. It follows that $[EF_{\max}^{(e)}:E]<\infty$. 
Hence there exists a finite separable extension $M/K$ such that
$F_{\max}^{(e)}\subset EM$, and this extension $M$ has the desired property. 
\end{proof}

\begin{proof}[Proof of Corollary \ref{hab2}] There exists a number field $k$ such that $G_A(k^\dagger K)$ is topologically finitely generated
by  Theorem \ref{main-thm}(d). By Proposition \ref{tocor} applied with $E=k^\dagger K$ there exists a finite extension $M/K$ such that 
$F\subset k^\dagger M$ for all $F\in \EEE_e(K_\tors/K)$. Replacing $M$ by $kM$ we get that $k^\dagger M\subset M_\ab$. 
\end{proof}

\begin{proof}[Proof of Corollary \ref{hab1}] From Theore \ref{main-thm}(b) and (c) we conclude that $G_A(K)$ is topologically finitely 
generated in the situation under consideration. The assertion now follows by Proposition \ref{tocor} applied with $E=K$. 
\end{proof}

\section{Application to the weak Hilbert property}\label{app}
In this section we use exactly the same notation as the paper \cite{bfp1} and the manuscript \cite{bfp2}. 
In particular, if $X$ is a normal variety over a field $K$ of characteristic zero and $W\subset X(K)$, then $W$ is said to be {\sl strongly thin}
if there exists a finite
family $(f_j: Y_j\to X)_{j=1,\cdots, s}$ of finite ramified morphisms and a proper closed subset $C$ of $X$ such that each $Y_j$ is normal and connected and such that 
$W\subset  C(K)\cup \bigcup_{j=1}^s f_j(Y_j(K))$. Furthermore $X$ is said to have {\sl the weak Hilbert property} (WHP for short), if $X(K)$ is not strongly thin. 

We shall strenghten and give a new proof of 
\cite[Theorem 1.7]{bfp1}, avoiding the original arguments from \cite[Section 5]{bfp1}. 
The new proof is based on Theorem \ref{main-theorem}, \cite[Corollary 4.3]{bfp2} and \cite[Theorem 1.4]{bfp1},  where \cite[Theorem 1.4]{bfp1} in turn relies on \cite{cdjlz} and builds on techniques from \cite{dan}. 

\begin{thm} \label{whp} Let $K$ be a number field and $A/K$ a geometrically simple abelian variety. Assume that $\rk(A(\Qq_\ab K))=\infty.$ Let $B/K$ be an abelian variety. Then 
$A_{K(B_\tors)}$ has property WHP over $K(B_\tors)$.
\end{thm}

\begin{proof} There exists a number field $k$ such that, if we put $E=kK$, the profinite group $G_B(E_{\ab})$ is topologically finitely generated 
(cf. Theorem \ref{main-thm} and Remark \ref{explO}). Let $K_\tors=K(B_\tors)$. 
Clearly $$G_B(E_{\ab})  = \Gal(K_\tors E_\ab/E_\ab)\cong  \Gal(K_\tors/E_\ab\cap K_\tors),$$
where the first equality follows by Remark \ref{expl}, so this group is topologically finitely generated. 
Let $F=E_\ab\cap K_\tors$. Now  $FE/E$ is an abelian Galois extension,  hence $K':=E\cap F$ is a finite extension of $K$ 
such that $\Gal(F/K')$ is abelian. 
Furthermore, due to the existence of the Weil pairing on $B$, the field $K_\tors$ contains all roots of unity. 
It follows that $\Qq_\ab K \subset F$. Thus $\rk(A(F))=\infty$. The following diagram shows the fields constructed so far.
$$\xymatrix{
     & E_\ab\ar@{-}[r]\ar@{-}[d] & E_\ab K_\tors\ar@{-}[dd]\\
E\ar@{-}[ur]\ar@{-}[d]\ar@{-}[r]   & FE\ar@{-}[d] \\
K'\ar@{-}[d] \ar@{-}[r]&F\ar@{-}[r] \ar@{-}[d]        & K_\tors\\
K\ar@{-}[r] & \Qq_\ab K\ar@{-}[ur] \\
}$$

By  \cite[Theorem 1.4]{bfp1} it follows that $A_F=(A_{K'})_F$ has WHP. From the finite generation of $\Gal(K_\tors/F)$ and \cite[Corollary 4.3]{bfp2}
we obtain the assertion. 
\end{proof} 

\begin{coro} (cf. \cite[Theorem 1.7]{bfp1})
Let  $A/\Qq$ be an elliptic curve. Let $B/\Qq$ be an abelian variety. Then 
$A_{\Qq(B_\tors)}$ has property WHP over $\Qq(B_\tors)$.
\end{coro}

\begin{proof} In that case $A$ is obviously geometrically simple. Furthermore we have $\rk(A(\Qq_\ab))=\infty$ by \cite[Lemma 2.1, Theorem 2.2]{fj}. The assertion now follows by 
Theorem \ref{whp}. 
\end{proof}

\begin{rema}\label{rem:zan}
Let us briefly explain the relation with the conjecture of Zannier in \cite[Section 2]{zan}. Consider the situation of Theorem \ref{whp} with $A=B$. Let $T=A(\ol{K})_\tors$.
Let $f: X\to A$ be 
a cover and assume $X$ is geometrically irreducible. Assume that $X$ is not isomorphic to an abelian variety.  Then $f$ is necessarily ramified and $M=f(X(K_\tors))$ is a
strongly thin set. By
Theorem \ref{whp} the complement of $M$ in $A(K_\tors)$ is Zariski dense. In particular the complement of $M\cap T$ in $A(K_\tors)$ is Zariski dense. 
The conjecture would predict the much stronger statement that $M\cap T$ is contained in a proper closed subset, and this without assumptions on the rank. Nevertheless 
one can see Theorem \ref{whp} as a result providing at least some evidence for the conjecture of Zannier. 
\end{rema}

\newpage

\small
{\sc Wojciech Gajda\\
Faculty of Mathematics and Computer Science\\
Adam Mickiewicz University\\
Uniwersytetu Pozna\' nskiego 4\\
61614 Pozna\'{n}, Poland}\\
E-mail adress: \texttt{\small gajda@amu.edu.pl}
\par\bigskip

{\sc Sebastian Petersen\\ 
Universit\"at Kassel\\
Fachbereich 10\\
Wilhelmsh\"oher Allee 71--73\\
34121 Kassel, Germany}\\
E-mail address: \texttt{petersen@mathematik.uni-kassel.de}

\end{document}